\documentclass[a4paper, 11pt]{amsart}
\usepackage{graphicx}

\begin{document}

\title{Quantum computing, Seifert surfaces and singular fibers}

\author{Michel Planat$\dag$, Raymond Aschheim$\ddag$,\\ Marcelo M. Amaral$\ddag$ and Klee Irwin$\ddag$}

\address{$\dag$ Universit\'e de Bourgogne/Franche-Comt\'e, Institut FEMTO-ST CNRS UMR 6174, 15 B Avenue des Montboucons, F-25044 Besan\c con, France.}
\email{michel.planat@femto-st.fr}

\address{$\ddag$ Quantum Gravity Research, Los Angeles, CA 90290, USA}
\email{raymond@QuantumGravityResearch.org}
\email{Klee@quantumgravityresearch.org}
\email{Marcelo@quantumgravityresearch.org}

\begin{abstract}

The fundamental group $\pi_1(L)$ of a knot or link $L$ may be used to generate magic states appropriate for performing universal quantum computation and simultaneously for retrieving complete information about the processed quantum states. In this paper, one defines braids whose closure is the $L$ of such a quantum computer model and computes their Seifert surfaces and the corresponding Alexander polynomial.
 In particular, some $d$-fold coverings of the trefoil knot, with $d=3$, $4$, $6$ or $12$,  define appropriate links $L$ and the latter two cases connect to the Dynkin diagrams of $E_6$ and $D_4$, respectively. In this new context, one finds that this correspondence continues with the Kodaira's classification of elliptic singular fibers. The Seifert fibered toroidal manifold $\Sigma'$, at the boundary of the singular fiber $\tilde {E_8}$, allows possible models of quantum computing.

 \end{abstract}

\maketitle

\vspace*{-.5cm}

\normalsize

\section{Introduction}

To acquire computational advantage over a classical circuit, a quantum circuit needs a non-stabilizer quantum operation for preparing a non-Pauli eigenstate, often called a magic state. The work about qubit magic state distillation \cite{Bravyi2004} was generalized to qudits \cite{Veitch2014} and multi-qubits (see \cite{Seddon2019} for a review). Thanks to these methods, universal quantum computation (UQC), the ability to prepare every quantum gate, is possible. A new approach of UQC, based on permutation gates and simultaneously minimal informationally complete positive operator valued measures (MICs), was  worked out in \cite{PlanatRukhsan,PlanatGedik}. It is notable that the structure of the modular group $\Gamma$ is sufficient for getting most permutation-based magic states \cite{PlanatModular} useful for UQC and that this can be thought of in terms of the trefoil knot $3$-manifold \cite{MPGQR1}.

It is desirable that the UQC approach of \cite{PlanatRukhsan}-\cite{MPGQR1} be formulated in terms of braid theory to allow a physical implementation. Braids of the anyon type, that are two-dimensional quasiparticles with world lines creating space-time braids, are nowadays very popular \cite{Freedman2003,Nayak2008}. Close to this view of topological quantum computation (TQC) based on anyons, we propose  a TQC based on the Seifert surfaces defined over a link $L$. The links in question will be those able to generate magic states appropriate for performing permutation-based UQC. 

In our previous work \cite{MPGQR1}, we investigated the trefoil knot $T_1=3_1$ as a possible source of $d$-dimensional UQC models through its subgroups of index $d$ (corresponding to $d$-fold coverings of the $T_1$ $3$-manifold) (see \cite[Table 1]{MPGQR1}). More precisely, the link $L7n1$, corresponding to the congruence subgroup $\Gamma_0(2)$ of the modular group $\Gamma$, builds a relevant qutrit magic state for UQC and is related to the Hesse configuration. The link $L6a3$, corresponding to the congruence subgroup $\Gamma_0(3)$ of $\Gamma$, builds a relevant two-qubit magic state and is related to the figure $GQ(2,2)$ of two-qubit commutation of Pauli operators. Then the link identified by the software SnapPy as $L6n1$ (or sometimes $L8n3$), corresponding to the congruence subgroup $\Gamma(2)$ of $\Gamma$, defines a $6$-dit MIC with the figure of Borromean rings as a basic geometry \cite[Fig.~4]{PlanatModular}. As shown below, it turns out that none of the two aforementioned links $L6n1$ and $L8n3$ are correctly associated to the subgroup $\Gamma(2)$ of $\Gamma$, but the link $6_3^3$ (related to the Dynkin diagram of $\tilde{D}_4$) is. The possible confusion lies in the fact that all three links share the same link group $\pi_1(L)$. Finally, the Dynkin diagram of $D_4$ (with the icosahedral symmetry of $H_3$ in the induced permutations) is associated to a $12$-dimensional (two-qubit/qutrit) MIC corresponding to the congruence subgroup $10A^1$ of $\Gamma$ \cite[Table 1]{PlanatModular}.

As announced in the abstract, we introduce a Seifert surface methodology for converting the UQC models based on the just described links into the appropriate braid representation permitted by Alexander's theorem. These calculations are described in Sec.~2. Then, in Sec.~3, we generalize the building of UQC models to affine Dynkin diagrams of type $\tilde{D_4}$, $\tilde{E_6}$ and $\tilde{E_8}$, that are singular fibers of minimal elliptic surfaces. Along the way, topological objects such as the $3$-torus, the Poincar\'e dodecahedral space \cite{Weeks2001} as well as the first amphicosm \cite{Chelnokov2017} are encountered. They are the precursors of $4$-manifold topology that is currently under active scrutiny \cite{Gompf1999}. Its possible role in models of UQC is discussed in the conclusion.

\section{Seifert surfaces and braids from $d$-fold coverings of the trefoil knot manifold (or of hyperbolic 3-manifolds)}

Alexander's theorem states that every knot or link can be represented as a closed braid \cite{AdamsBook}. A Seifert surface $F$ of a knot $K$ or a link $L$ is an oriented surface within the $3$-sphere $S^3$ whose boundary $\partial F$ coincides with that knot or link. Given a basis $\{f_k\}$ for the first homology group $H_1(F:\mathbb{Z})$ of $F$, one defines a Seifert matrix $V$ whose $(i,j)$-th entry is the linking number of the component $f_i$ and the positive push-off $f_j^+$ of the component $f_j$ along a vector field normal to $F$. Then a useful invariant of $L$ is the (symmetrized) Alexander polynomial \cite{AdamsBook},\cite[Sec.~2.7]{Akbulut2016}  

\begin{equation}
\Delta_L(t)=t^{-r/2} \det(V-tV^T),
\label{eq1}
\end{equation}
with $V^T$ the transpose of $V$ and $r$ the first Betti number $F$. By definition $\Delta_L(t^{-1})=\Delta_L(t)$.

\begin{figure}[ht]
\centering 
\includegraphics[width=6cm]{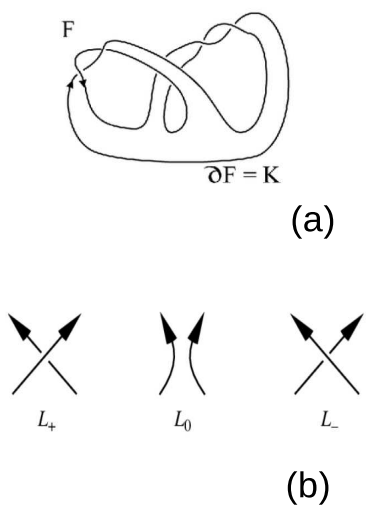}
\caption{(a) The Seifert surface $F$ for a trefoil knot $K$, (b) the types of crossings for the skein relation of a link $L$.}
\label{skein}
\end{figure}

There exists a remarkable topological property  of $\Delta_L(t)$ called a skein relation \footnote{The Jones polynomial used for defining anyons obeys a different skein relation than the Alexander polynomial \cite{Kauffman1987} so that the rules for braiding are also different from those resulting from the Seifert surfaces.}. If $L_+$, $L_0$ and $L_-$ are links in $S^3$, with projections differing from each other by a single crossing, as in Fig.~\ref{skein}b, then 

\begin{equation}
\Delta_{L_+}(t)-\Delta_{L_-}(t)=(t^{1/2}-t^{-1/2}) \Delta_{L_0}(t).
\label{eq2}
\end{equation}

When $L$ is a knot $K$, there is a connection of $\Delta_K(t)$ with a combinatorial invariant $\nu$ of the $3$-manifold $S_K^3$ obtained from the $0$-surgery 
along $K$ in $S^3$ as follows

\begin{equation}
\nu(S_K^3)=\frac{\Delta_K(t)}{(t^{1/2}-t^{-1/2})^2}.
\label{eq3}
\end{equation}

The invariant $\nu(S_K^3)$ is called Milnor (or Reidemeister) torsion \cite{Milnor1962}. The main interest of $\nu$ is its ability to distinguish closed manifolds which are homotopy equivalent while being non homeomorphic.

 The Seifert surface can be drawn  from the braid representation. A good reference is \cite{Wijk2006} and the related website \cite{JuliaCollins}. The software SeifertView provides a visualization of the Seifert surface. Below, to practically obtain the braid representation and the corresponding Alexander polynomial, we proceed as follows. 
With the software SnapPy \cite{SnapPy}, one defines the link from its name [e.g. {\it M=Manifold}(\lq $3_1$') for the trefoil knot $K=T_1$] or from its PD representation available after drawing the link in the pink editor [e.g. {\it trefoil=}$[(6,4,1,3),(4,2,5,1),(2,6,3,5)]$, {\it L=Link(trefoil)} and
 {\it L.braid\_word()} for obtaining the braid associated to $T_1$ as $[-1, -1, -1]$ and {\it L.braid}\_{\it matrix()} for obtaining the Seifert matrix $V$]. Then, with Magma software \cite{Magma}, the Alexander polynomial follows [{\it det:=Determinant} $(u*V-v*Transpose(V)$){\it; det;} to obtain $u^2 - u*v + v^2$, or $t-1+t'$ after replacing $u$ by $t^{1/2}$ and $v$ by $t'=t^{1/2}$].

Results are summarized in Table 1.

\scriptsize
\begin{table}[h]
\begin{center}
\begin{tabular}{|c|l|c|c|r|}
\hline 
\hline
 source  &  target & MIC & braid word& Alexander polynomial  \\
\hline
trefoil & L7n1 & QT Hesse  & $(ab)^3b$ & $t^{5/2}-t^{3/2}+t'^{(3/2)}-t'^{(5/2)}$ \\
.       & L6a3 & 2QB Doily &  $ABCDCbaCdEdCBCDCeb$ &$-3t^{1/2}+3t'^{(1/2)}$\\
.           & $6_3^3$ & 6-dit-MIC & $(ab)^3$ & $t^2-t-t'+t^2$\\
.       &\tiny{$D_4$ Dynkin} & \tiny{2QB-QT MIC} & $ABCCbaCCBCCb$ & \tiny{ $-t^{3/2}+3t^{1/2}+t'^{(3/2)}-3t'^{(1/2)}$}\\
\hline
\hline
fig. eight &\tiny{L14n55217}&  7-dit MIC &   $AbbcbcbDacBacdcb$  & $-t^4+7t^3-11t^2+8t-6+\cdots $\\
\hline
Whitehead & \tiny{L12n1741}&QT Hesse& \tiny{AbcDEFeDCBDacBdcdEdfCbdCddddeD}&$-2t^3+6t^2-6t+4+ \cdots$\\ 
. & \tiny{L13n11257} & 5-dit MIC & AbCCbDaCBcDcDcbCD& $t^{9/2}-6t^{7/2}+15t^{5/2}-21t^{3/2}$\\
 &  &  & &  $+21t^{1/2}+\cdots$\\
\hline
$6_3^2=L6a1$ & \tiny{L12n2181} & QT Hesse & \tiny{ABcdEFceGbdFaedCBcdEdfcEgbdfedc}& $4t^{5/2}-12t^{3/2}+16t^{1/2}+\cdots$  \\
. & L14n63905 & 2QB MIC & \tiny{AbCddEdFedcBdaEdfCbceDccDcBC}& $t^4-7t^3+22t^2-41t+50+\cdots$\\
\hline
magic L6a5 & \tiny{L14n63788} & QT Hesse & ABCdEEEFEDcebdacEbEED  & $t^4-2t^3+2t-2+\cdots$\\
 &  &  &  ceDefedCeBdCEDe &  \\
\hline
\hline
\end{tabular}
\vspace*{.5cm}
\caption{A few models of universal quantum computation (UQC)  \cite{MPGQR1,MPGQR2} translated into the language of braids and their Seifert surfaces. The source is a knot (such as the trefoil knot) or a link and the target is a link $L$ associated to a degree $d$ covering of $L$-manifold that defines an appropriate magic state for UQC and a corresponding MIC (minimal informationally complete) POVM. Cases $d=3,4,5,...$ correspond to the geometry of the Hesse configuration, to the doily $GQ(2,2)$ finite geometry, to the Petersen graph \ldots The notation for the braids is that of \cite{Wijk2006}. the notation $t'$ means $t^{-1}$.}
\end{center}
\end{table}

\normalsize

Before going further, let us recall the homomorphism between the conjugacy classes of subgroups of index $d$ of a group $G$ and the $d$-fold coverings of a manifold $M$ whose fundamental group is $G=\pi_1(M)$ \cite{Mednykh2006}. This relationship is not one-to-one (not an isomorphism) in the sense that a $\pi_1(M)$ may characterize distinct manifolds $M$. A simple example is for the $d$-coverings of manifold with characteristic $\nu=2g-2=0$ (with $g$ the genus) whose number is the sum of divisor function  $\sigma(d)$ \cite[Sec.~3.4]{LiskovetsMednykh2009}. The same cardinality structure of finite subgroups of $\pi_1(L)$ [denoted $\eta_d(L)$ in our previous paper \cite{MPGQR1}] is also encountered for the distinct links $L6n1$, $L8n3$, $6_3^3$ and for the Kirby link introduced later. 

\subsection{The braids built from the trefoil knot that are associated to the qutrit link L7n1 and the two-qubit link L6a3}

\noindent

\begin{figure}[ht]
\centering 
\includegraphics[width=6cm]{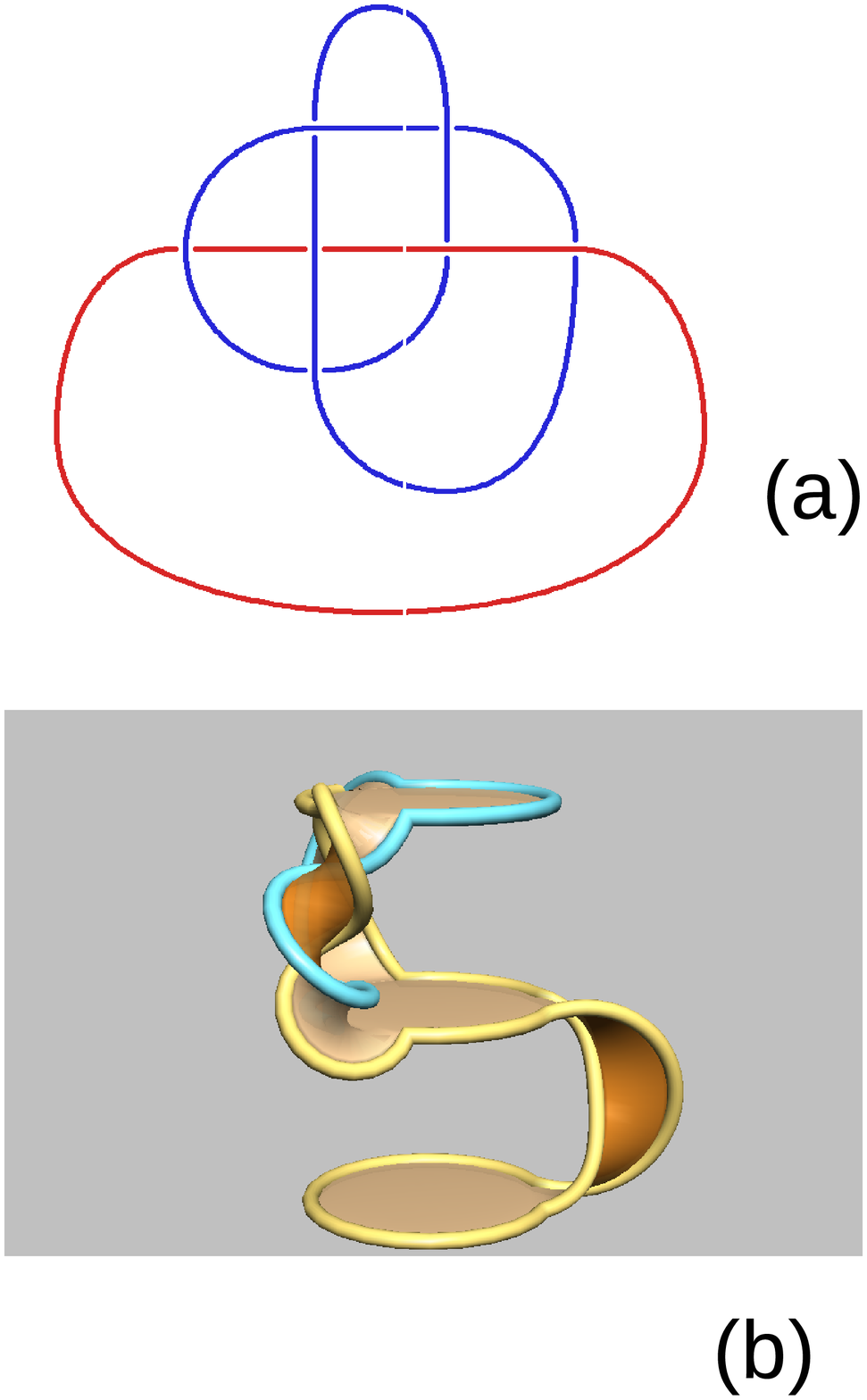}
\includegraphics[width=6cm]{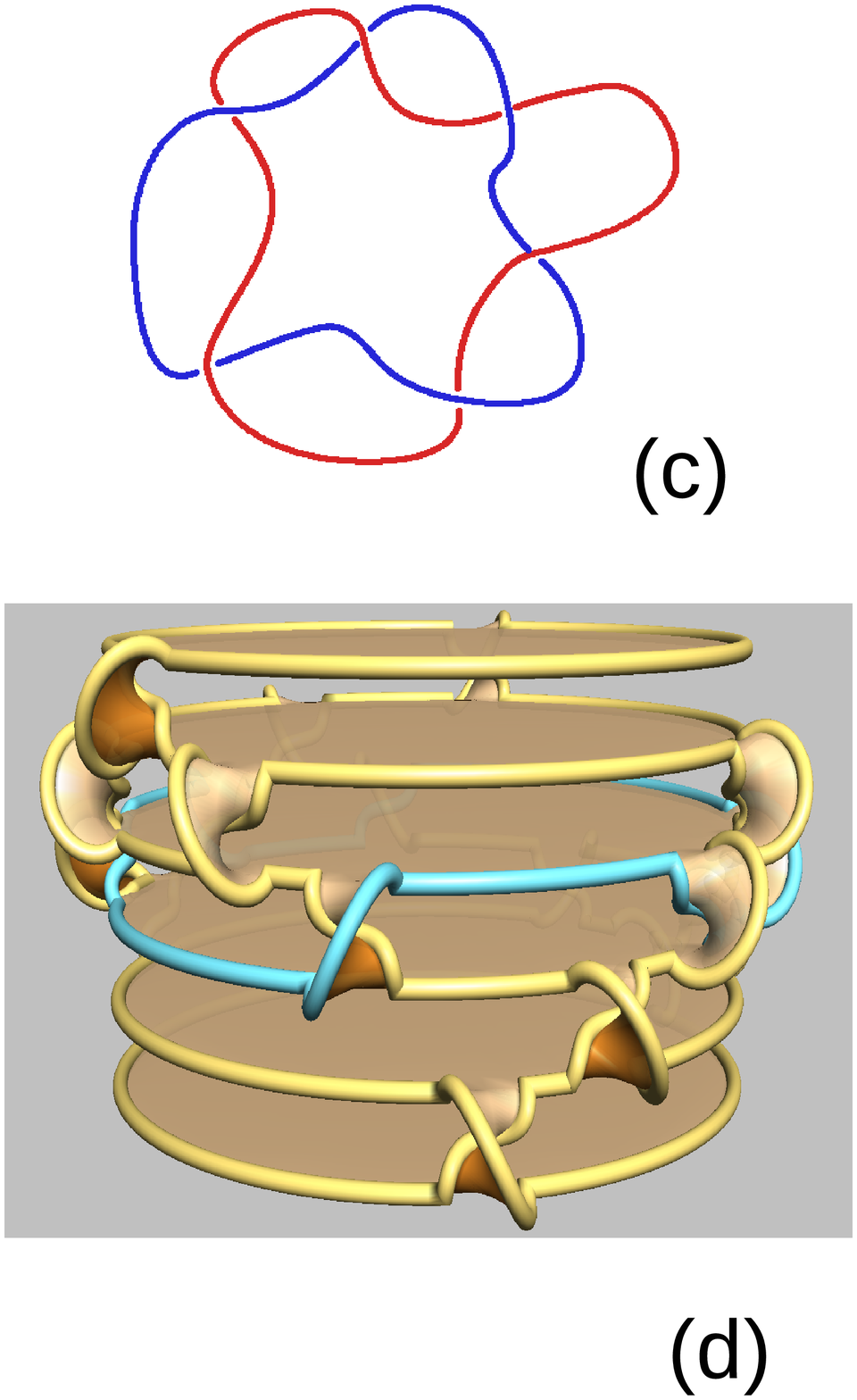}
\caption{(a) The link L7n1 defining the qutrit MIC and (b) its Seifert surface. (c) The link L6a3 defining the two-qubit MIC and (d) its Seifert surface.}
\label{fig1}
\end{figure}

As announced in the introduction, one refers to \cite[Table 1]{MPGQR1} that lists the topological properties of $d$-fold coverings $d=1 \ldots8 $ of trefoil knot manifold. as obtained from SnapPy and also identifies the corresponding congruence subgroups of $\Gamma$ previously investigated in \cite{PlanatModular}. 
From now, one denotes $A=a^{-1}, B=b^{-1}$,\ldots and (.,.) means the group theoretical commutator of the entries.
The link $L7n1$ corresponds to the congruence subgroup $\Gamma_0(2)$ of $\Gamma$, its fundamental group $\pi_1=\left\langle a,b|(a,B^2)\right\rangle$ builds a qutrit magic state for UQC of the type $(0,1,\pm 1)$ and a MIC with the Hesse geometry \cite[Fig.~1a]{MPGQR1}. Fig.~\ref{fig1}a is the drawing of $L7n1$ and Fig.~\ref{fig1}b is that of the Seifert surface for the braid word $(ab)^3b$.

The link $L6a3$ corresponds to the congruence subgroup $\Gamma_0(3)$ of $\Gamma$, its fundamental group $\pi_1=\left\langle a,b|(a,b^3)\right\rangle$ builds a two-qubit magic state for UQC of the type $(0,1,-\omega_6, \omega_6-1)$, $\omega_6=\exp(\frac{2i\pi}{6})$,  as well as a MIC with the geometry of the generalized quadrangle of order two $GQ(2,2)$ \cite[Fig.~1b]{MPGQR1}.  Fig.~\ref{fig1}c is the drawing of $L6a3$ and Fig.~\ref{fig1}d is that of the Seifert surface for the braid word $ABCDCbaCdEdCBCDCeb$.

The Alexander polynomials are made explicit in Table 1.

\subsection{The braid built from the trefoil knot that is associated to the $6$-dit link $6_3^3$ and related braids with the same fundamental group}

\noindent

There are eight conjugacy classes of subgroups of index $6$ of $\Gamma$ corresponding to eight $6$-fold coverings over the trefoil knot manifold. They are listed and identified in \cite[Table 1]{MPGQR1}. We are first interested in the unique regular covering $M$ of degree $6$ with homology $Z+Z+Z$ and three cusps corresponding to the congruence subgroup $\Gamma(2)$. The cardinality sequence of subgroups for the fundamental group of this particular covering is that of the link $6_3^3$

\begin{equation}
\eta_d(6_3^3)=[1, 7, 16, 60, 122, ~794, 4212, 35276, 314949, \ldots]
\label{eqn4}
\end{equation}

It turns out that every degree $6$ covering of the trefoil manifold leading to a $6$-dit MIC [with magic state of the type $(0,1,\omega_6-1,0,-\omega_6,0)$] share the same fundamental group. For $M$, SnapPy randomly provides several choices such as $L6n1$ or $L8n3$ that, of course, share the same cardinality sequence as $6_3^3$. In the Knot Atlas at \lq http://katlas.org/wiki/L6n1', one finds the sentence \lq $L6n1$ is $6^3_3$ in Rolfsen's table of links'. But that seems to be a wrong statement.

\begin{figure}[ht]
\centering 
\includegraphics[width=6cm]{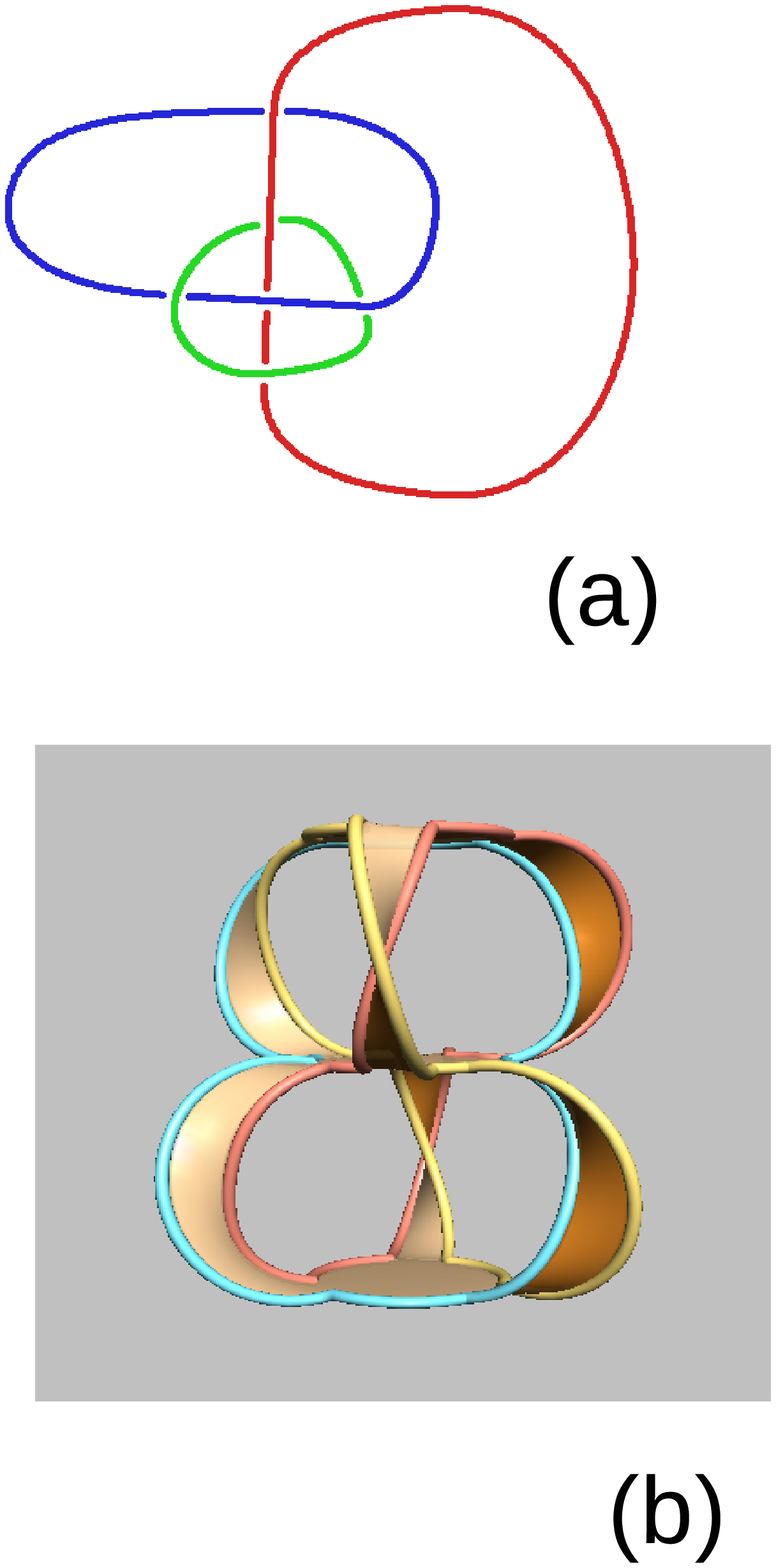}
\includegraphics[width=6cm]{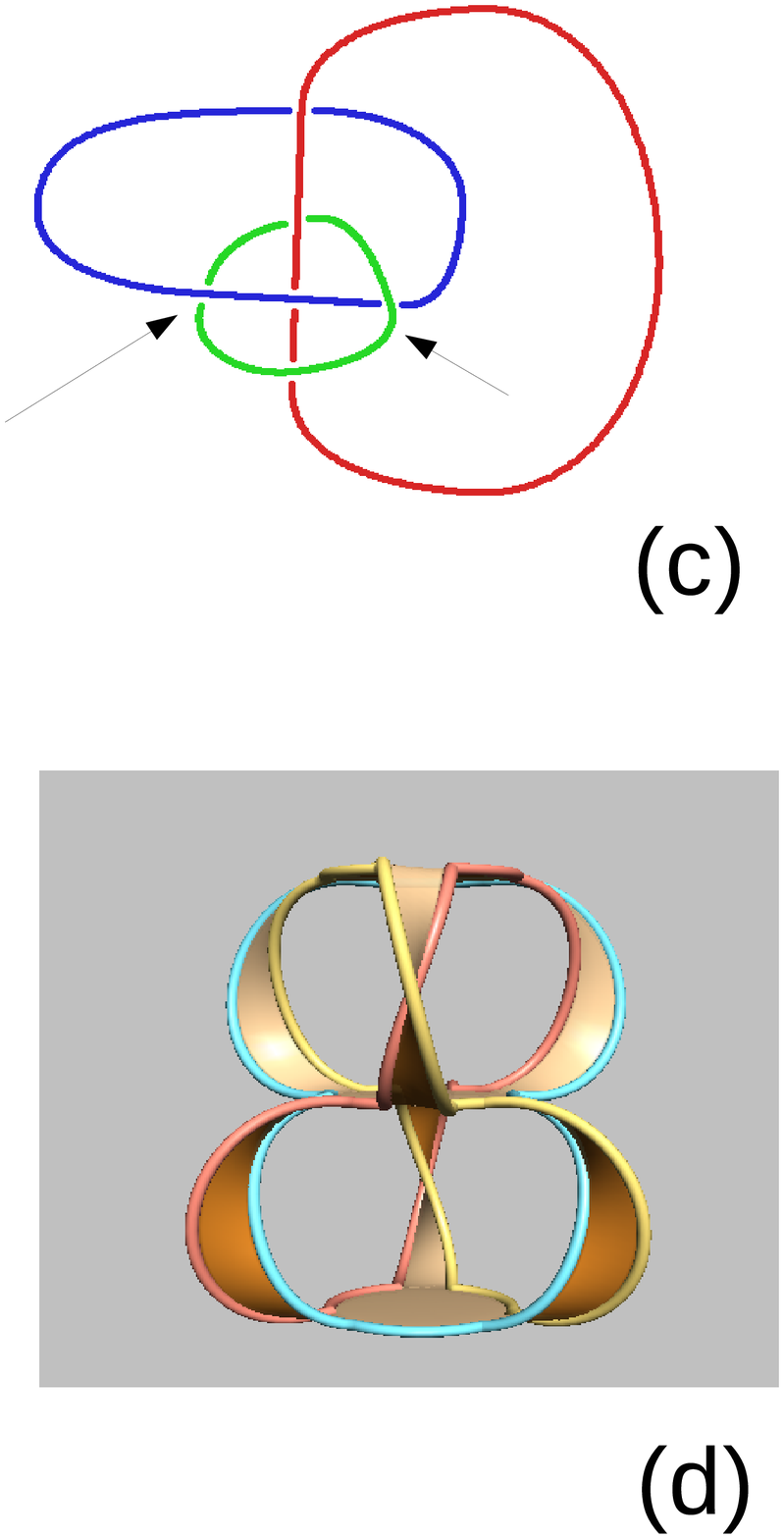}
\caption{(a) The link $6_3^3$ corresponding to the 6-dit MIC and the congruence subgroup $\Gamma(2)$ of $\Gamma$ and (b) its Seifert surface.
(c) The Kirby link $L_K$ (see the arrows for the up/down changes) and its Seifert surface (d) (observe the color changes).}
\label{fig2}
\end{figure}

Let us observe that by performing $(-2,1)$ surgery on all cusps of $M$ or of $6_3^3$ and introducing the Dynkin diagram $\tilde{E}_6$ of affine $E_6$ (see Fig.~\ref{fig4}b in Sec.~3 for details) one gets
\begin{eqnarray}
&\eta_d[M(-2,1)]=\eta_d[6_3^3(-2,1)]=\eta_d(\tilde{E}_6) \nonumber \\
&=[ 1, 1, 4, 2, 1,~ 6, 3, 2, 10, 1,~ 1, 19, 3, 3, 14,~ 3, 1, 36, 3, 2, \ldots ],\nonumber \\ \nonumber
\end{eqnarray}
while performing $(-2,1)$ surgery on cusps of $L6n1$ and introducing the Dynkin diagram of $E_6$ one gets 
\begin{eqnarray}
&\eta_d[L6n1(-2,1)]=\eta_d(E_6)=\eta_d(2T)\nonumber \\
&= [ 1, 0, 1, 1, 0, 1, 0, 1, 0, 0, 0, 1, 0, 0, 0, 0, 0, 0, 0, 0, 0, 0, 0, 1, 0, \ldots ],\nonumber \\ \nonumber
\end{eqnarray}
where $2T$ is the binary tetrahedral group.

We do not provide the result of performing $(-2,1)$ surgery on $L8n3$ which provides a still different result that we could not identify. 
We conclude that the correct identification of the manifold $M$ should be $6_3^3$ although we do not yet have a rigorous proof.

The Seifert surface for the link $6_3^3$ is drawn in Fig.~\ref{fig2}a and the corresponding Seifert surface associated to braid word $(ab)^3$ is Fig.~\ref{fig2}b.
Switching the up/down positions of circles at two points (as shown in Fig.~2c) provides the Kirby link $L_K$ drawn in \cite[Fig.~3]{Kirby1979} that we reproduce in  Fig.~\ref{fig2}c. [Applying the surgeries as $L_K(4,1)(1,1)(2,1)$ to red, blue and green circles, one gets the Brieskorn sphere $\Sigma(2,3,5)$, alias the Poincar\'e dodecahedral space]. The corresponding Seifert surface associated to the braid word $aBabAb$ is shown in Fig.~\ref{fig2}d. It is notable that $(-2,1)$ surgery on $K$ leads to the Dynkin diagram for $A_3$ of Weyl group $S_4$. Both links $6_3^3$ and $K$ are described by the same Alexander polynomial $t^2-t-t'+t^2$.

\subsubsection*{The $6$-cover of the trefoil knot manifold corresponding to the congruence subgroup $3C^0$ of $\Gamma$}

Let us conclude this subsection by another observation concerning the $6$-cover (that we denote  $M'$) of the trefoil knot manifold identified in \cite[Table 1]{MPGQR1} corresponding to the congruence subgroup $3C^0$. Again, the cardinality sequence of subgroups of $\pi_1(M')$ is that of $6_3^3$, L6n1, L8n3  or $L_K$ but $M'$ can be distinguished from the manifolds corresponding to these links since one gets under $-2$-surgery $\eta_d[M'(-2,1)]=\eta_d[\tilde{D}_4]$, where $\tilde{D}_4$ is the Dynkin diagram of affine $D_4$ (as well as the smallest elliptic singular fiber of Kodaira's classification, see Fig.~\ref{fig4}a).

\subsection{The braid built from the trefoil knot that is associated to the two-qubit/qutrit MIC with icosahedral symmetry of the permutation representation}

\noindent

\begin{figure}[h]
\centering 
\includegraphics[width=6cm]{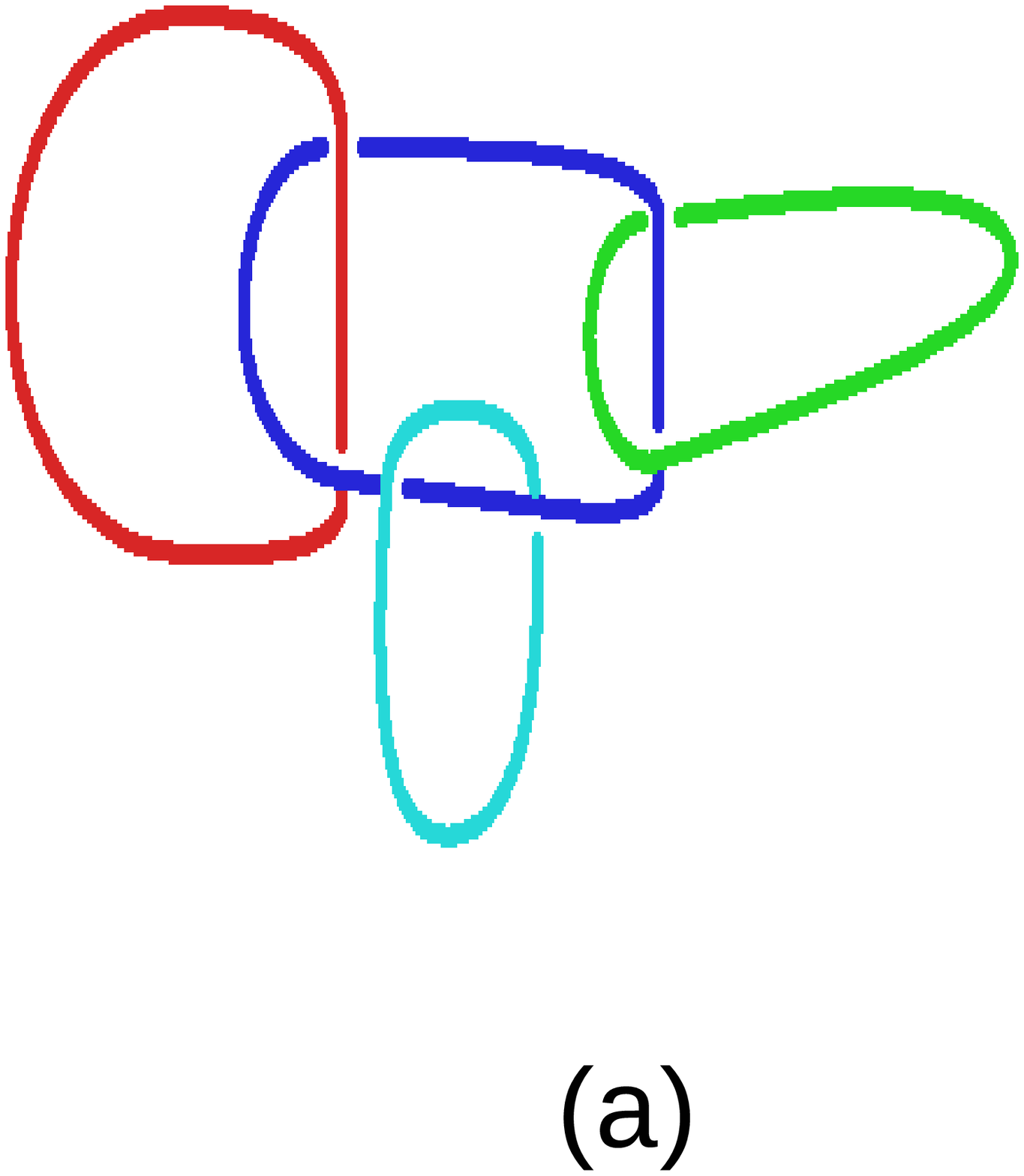}
\includegraphics[width=6cm]{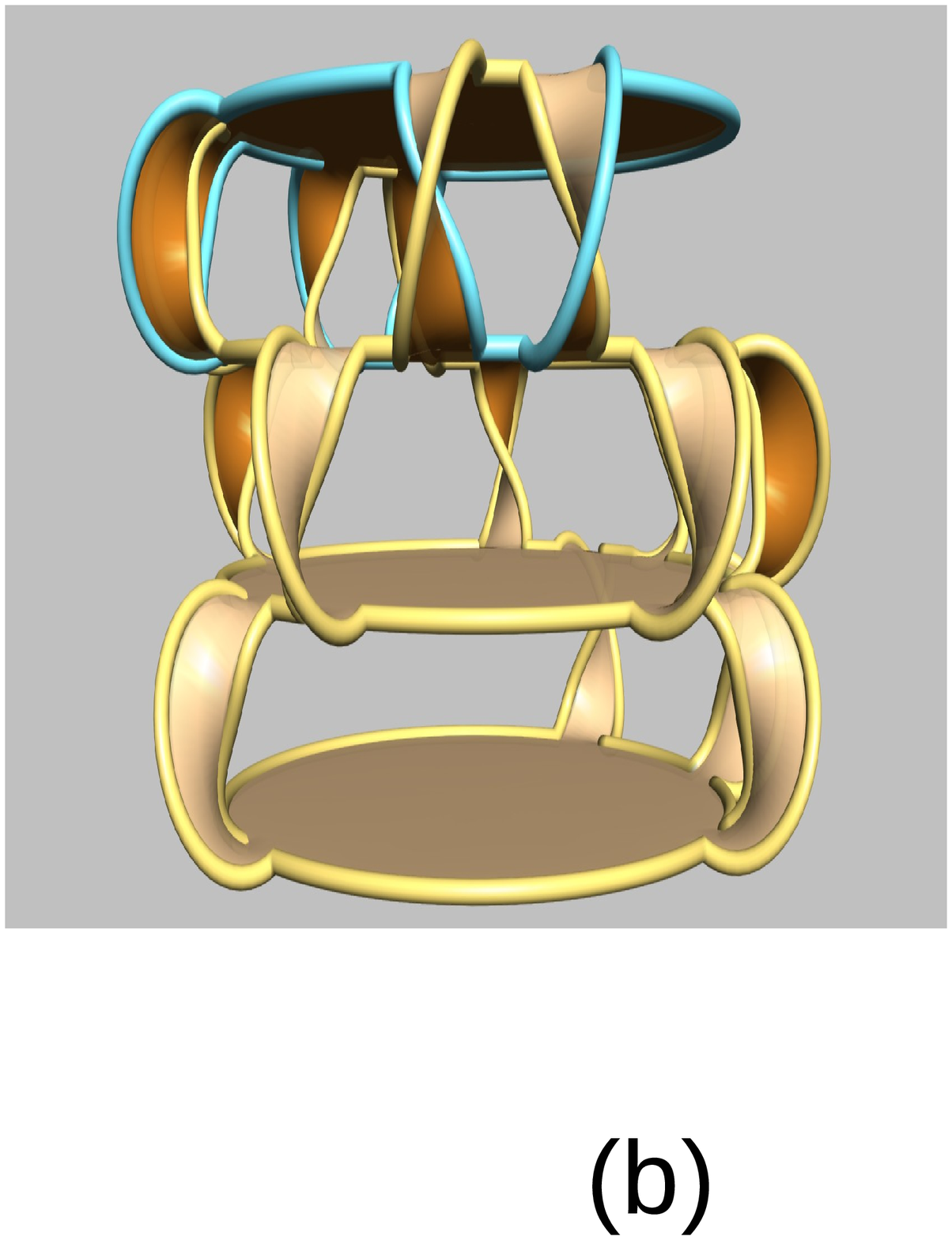}
\caption{(a) The Dynkin diagram for the $D_4$ manifold attached to the 2QB-QT MIC consists of pairs of handles arranged as the Coxeter diagram of $D_4$, (b)
the Seifert surface corresponding to the Dynkin diagram of $D_4$ and braid word $ABCCbaCCBCCb$.}
\label{fig3}
\end{figure} 

In \cite[Table 1]{PlanatModular}, the first author identified a two-qubit/qutrit MIC corresponding to the congruence subgroup $10A^1$ of the modular group $\Gamma$. The geometry of this MIC is that of the four-partite graph $K(3,3,3,3)$. More precisely , the subgroup of index $12$ corresponding to $10A^1$ is generated as $G=\left\langle a,b,c,d|(a,b),(a,c),(a,t)\right\rangle$ that describes the Dynkin diagram of $D_4$. It builds a magic state for UQC of the type $(1,1,0,0,0,0,-\omega_6,-\omega_6,0, \omega_6-1,0,\omega_6-1)$. Under $(-2,1)$ surgery on all four links one recovers the quaternion group.  The permutation representation that organizes the cosets of $10A^1$ in $\Gamma$ is the icosahedral group $\mathbb{Z}_2 \times PSL(2,5)$, alias the Coxeter group of the Dynkin diagram $H_3$. The braid word and the corresponding Alexander polynomial are given in Table 1, the Seifert surface is in Fig.~\ref{fig3}.

\subsection{Braids from $d$-fold coverings of hyperbolic $3$-manifolds}

\noindent

Models of UQC from MICs are also sometimes associated to links as already recognized in \cite{MPGQR1,MPGQR2}. Some of them are listed in Table 1 together as a corresponding braid word and Alexander polynomial. These models happen to be more complicated. We do not describe them in more detail.

\section{Quantum computing from affine Dynkin diagrams}

\begin{figure}[ht]
\centering 
\includegraphics[width=4cm]{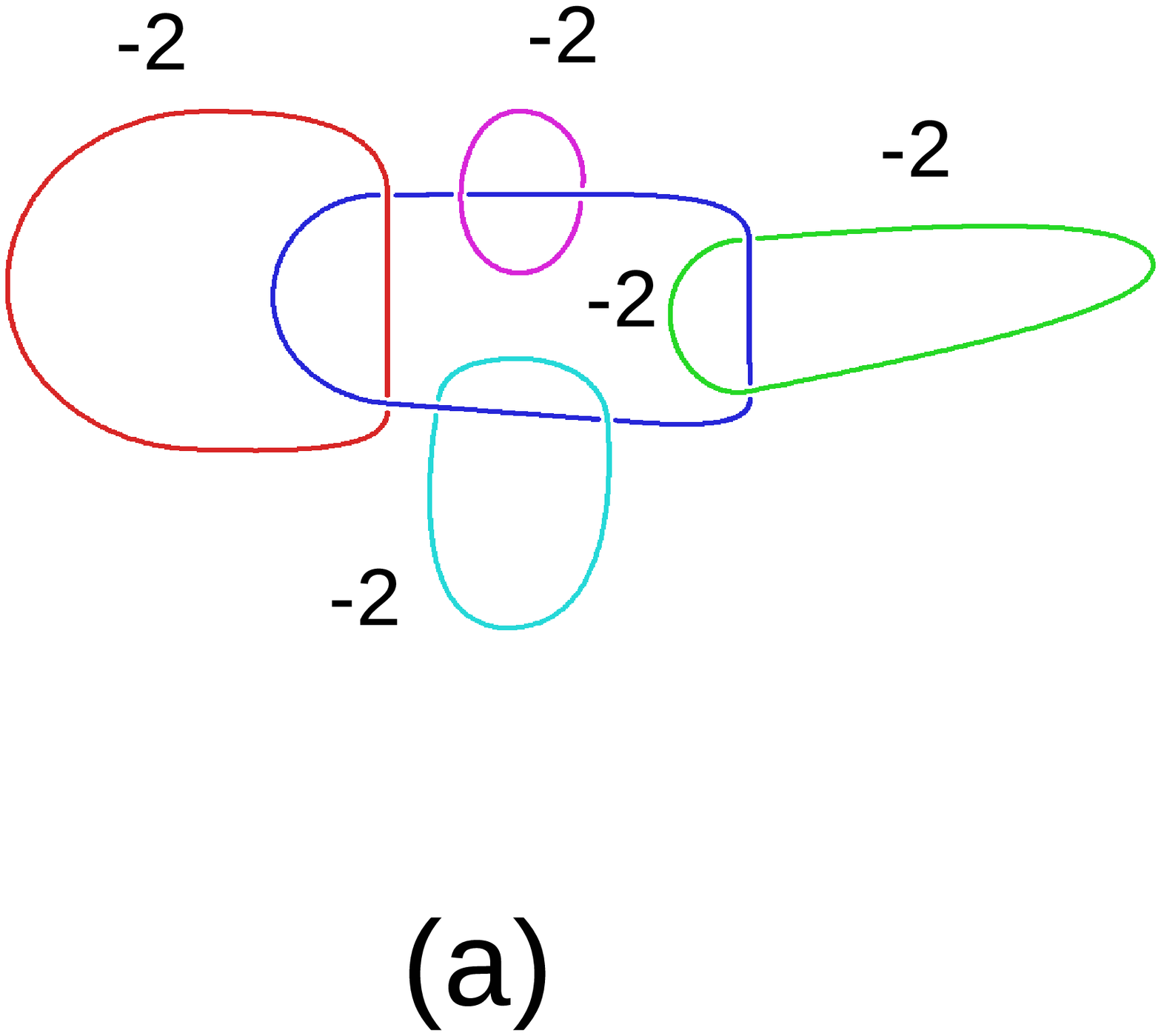}
\includegraphics[width=4cm]{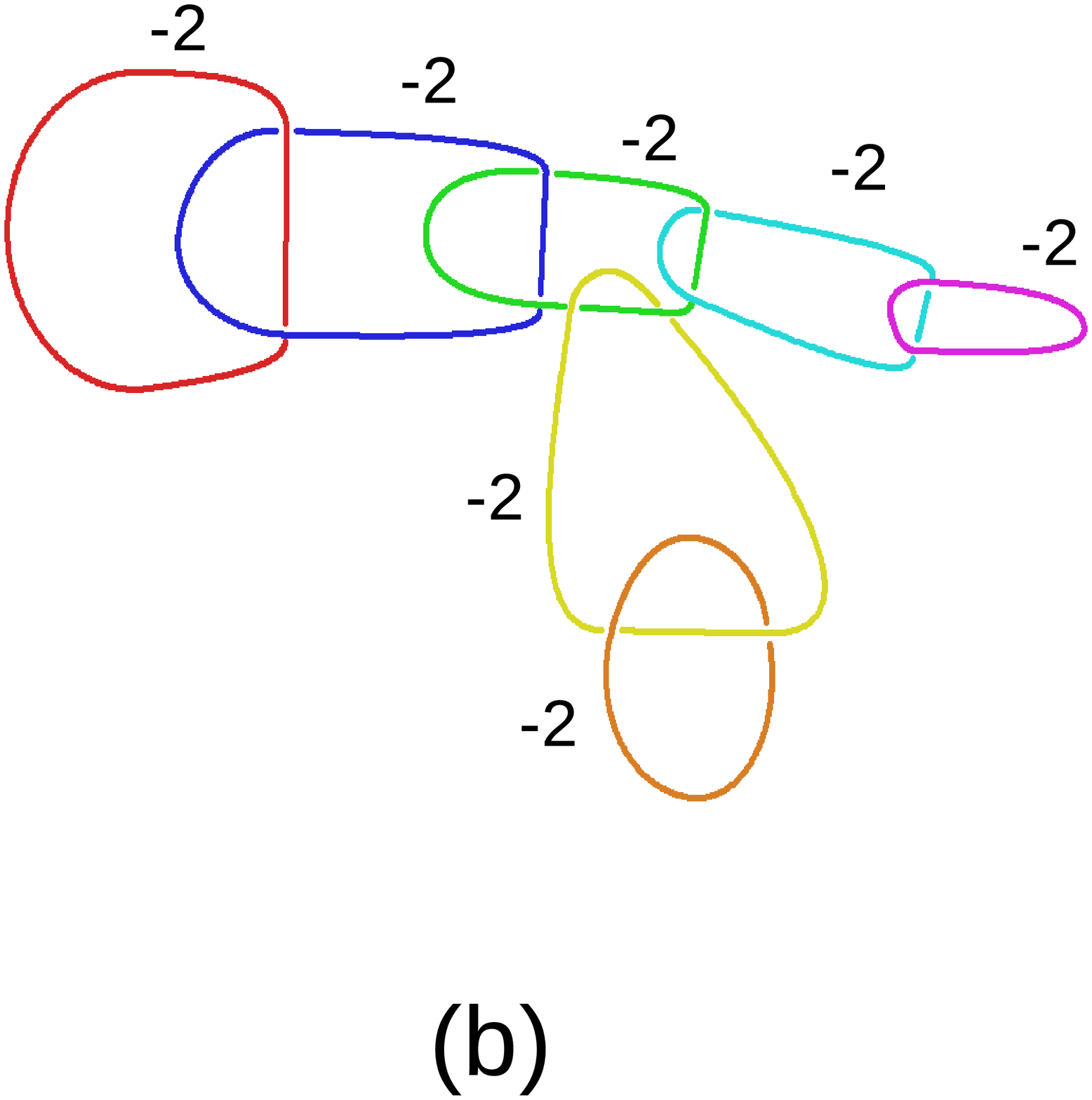}
\includegraphics[width=4cm]{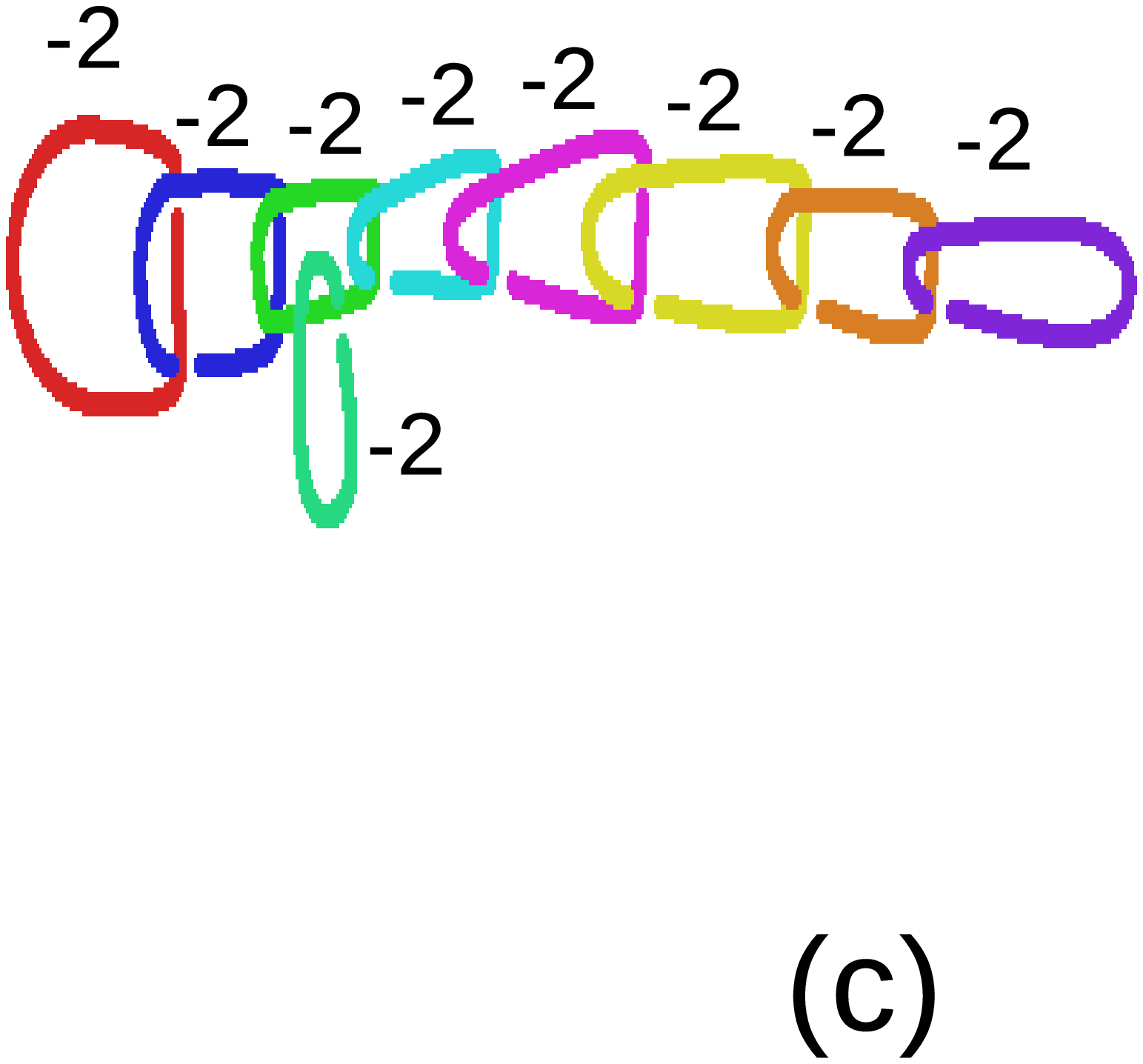}
\caption{ A few singular fibers in Kodaira's classification of minimal elliptic surfaces. (a) Fiber $I_0^*$ (alias $\tilde{D}_4$), (b) fiber $IV^*$ (alias $\tilde{E}_6$) and (c) fiber $II^*$ (alias $\tilde{E}_8$). }
\label{fig4}
\end{figure}

In the previous section, we found that some MICs for UQC (our approach of universal quantum computing with complete quantum information) relate to Coxeter-Dynkin diagrams: $\tilde{E_6}$ and $\tilde{D_4}$ for the $6$-dit MIC and $D_4$ for the two-qubit/qutrit MIC. This is an unexpected observation that we would like to complete by another one: the possibility of defining UQC from the singular fiber $II^*=\tilde{E_8}$ of Kodairas's classification of minimal elliptic surfaces (see Fig.~\ref{fig4}c). This classification is useful in the understanding of $4$-manifold topology as shown in \cite[p. 320]{Scorpian2005}, see also \cite{Aschheim2019} for a different perspective.
 
Taking $\pi_1(\tilde{E}_8)$  as the fundamental group of affine $E_8$, the subgroup structure of $\pi_1(\tilde{E}_8)$ has the following cardinality list

\begin{equation}
\eta_d(\tilde{E}_8)=[1,1,{\bf 2},{\bf 2},1,~{\bf 5},{\bf 3},2,4,1,~1,12,3,3,{\bf 4},\ldots]
\end{equation}

 where the bold characters mean that one of the subgroup leads to a MIC, as in \cite{MPGQR1}. It is worthwhile to observe that the boundary of the manifold associated to $\tilde{E}_8$ is the Seifert fibered toroidal manifold \cite{Wu2012}, denoted $\Sigma'$ in \cite[Table 5]{MPGQR1}. It may also be obtained by $0$-surgery on the trefoil knot $T_1$.

For the sequence above the coverings are
\begin{eqnarray}
&[\tilde{E}_8,\tilde{E}_6,\{\tilde{D}_4,\tilde{E}_8\},\{\tilde{E}_6,\tilde{E}_8\},\tilde{E}_8,~\{BR_0,\tilde{D}_4,\tilde{E}_6\},\{\tilde{E}_8\},\{\tilde{E}_6\},\{\tilde{D}_4,\tilde{E}_8\},\tilde{E}_6,\nonumber \\
&\tilde{E}_8,\{BR_0,\tilde{D}_4,\tilde{E}_6,\tilde{E}_8\},\{\tilde{E}_8\},\{\tilde{E}_6\},\{\tilde{D}_4,\tilde{E}_8\},\cdots] \nonumber \\ \nonumber
\end{eqnarray}
One observes that the subgroups/coverings are fundamental groups for  $\tilde{E}_8$ $\tilde{E}_6$, $\tilde{D}_4$ or $BR_0$,
where $BR_0$ is the manifold obtained by $0$-surgery on all circles of Borromean rings. The cardinality sequence of subgroups of $BR_0$ is 

\begin{equation}
\eta_d(BR_0)=[1, 7, 13, 35, 31, 91, 57, 155, 130, 217,\cdots]
\label{eqBR0}
\end{equation}

which is recognized as A001001 in Sloane's encyclopedia of integer sequences with the title \lq Number of sublattices of index $d$ in generic $3$-dimensional lattice'.

As already given in Sec.~1, the subgroup structure of $\pi_1(\tilde{E}_6)$ has the following cardinality list

\begin{equation}
\eta_d(\tilde{E}_6)=[ 1, 1, 4, {\bf 2}, 1,~ 6, 3, 2, 10, 1,\ldots].
\end{equation}

For this sequence the coverings are
$$[\tilde{E}_6,\tilde{E}_6,\{BR_0,\tilde{E}_6\},\tilde{E}_6,\tilde{E}_6,~\{BR_0,\tilde{E}_6\},\{\tilde{E}_6\},\{\tilde{E}_6\},
\{BR_0,\tilde{E}_6\},\tilde{E}_6 \cdots].$$

The subgroup structure for $\tilde{D}_4$ is

\begin{equation}
\eta_d(\tilde{D}_4)=[ 1, 7, 5, 23, 7,~ 39, 9, 65, 18, 61, \cdots]
\end{equation}

which corresponds to A263825 in Sloane's encyclopedia of integer sequences with the title 
\lq Total number $c_{\pi_1(B_1)}(n)$ of n-coverings over the first amphicosm $B_1$' \cite{Chelnokov2017}. 

To be exhaustive, let us mention that $\tilde{E}_7$ is $III^*$ Kodaira's singular fiber.
Following \cite[Table 1]{Kirby1999}, it can be obtained by $(-2,1)$-surgery on the link $L4a1$ \footnote{In \cite[Table 1]{Kirby1999}, $II^*$ is $3_1(0,1)$, $III^*$ is $L4a1(-2,1)$ and $IV^*=6_3^3(-2,1)$, as one expects.}. Observe that $L4a1$ has the same fundamental group as $L7n1$ (see Sec.~2.1) and $\eta_d[L4a1(-2,1)]=\eta_d(\tilde{E}_7)=[1,3,1,7,3,~5,1,16,2,11,\cdots]$.

$\tilde{E}_7$ has coverings of type $BR_0$, $\tilde{D}_4$ and $\tilde{E}_7$.

\subsubsection*{Reidemeister torsion of $\Sigma'$}

\noindent

As a final note for this section, according to (\ref{eq3}), the \lq twisted' Reidemeister torsion $\nu_t$ of the $3$-manifold obtained from $0$-surgery
 along a knot $K$ in $S^3$ is the Alexander polynomial of $K$. Thus for the trefoil knot $T_1=3_1$ one gets $ \nu_t(S^3_{T_1})=\nu_t(\Sigma')=t-1+t'$.

Let $BR=L6a4$ be the Borromean rings, the manifold $BR_0$ as above (obtained by $0$-surgery on all circles of $BR$) and $BR_1$ be the manifold obtained by $0$-surgery on two circles of $BR$. The cardinality sequence $\eta_d(BR_0)$ is as in (\ref{eqBR0}) and $\eta_d(BR_1)$ is found as in (\ref{eqn4}). In principle, one can compute the Reidemeister torsion for $BR_0$ and $BR_1$ \cite[Sec 2.4]{Nicolaescu}. This is left open in this paper.

\section{Conclusion}

In the first part of this paper, it has been found that some coverings of the trefoil knot, or of other knots or links useful for informationally complete UQC, are $3$-manifolds originating from a link within the $3$-sphere. Such links have Seifert surfaces and well defined braids. We pointed out some coverings of the trefoil knot (of index $6$ and $12$, respectively) connecting to the symmetry of Dynkin diagrams ( affine $E_6$ and $D_4)$.

In the second part, it has been found that the  singular fiber $\tilde{E}_8$ may be used to generate UQC and that its coverings are the singular fibers $\tilde{E}_6$ and $\tilde{D}_4$ (alias the first amphicosm), as well as the manifold $BR_0$ obtained by $0$-surgery on all circles of Borromean rings. Generalizing the analysis to coverings of $\tilde{D}_4$ and $BR_0$, one finds that the $3$-torus enters the game \cite{Weeks2001}. The full understanding of these facts needs the theory of $4$-manifolds \cite{Kirby1979,Scorpian2005}.

There exists surface braids wrapped around singularities within $4$-manifolds that support quasiparticles closely related to anyon states \cite{Atiyah2017}. More work is necessary to connect our observations to this latter work. Another connection seems to be the theory of quasicrystals that uses $E_8$ symmetry, e.g. \cite{FangIrwin2016, GQRgroup2016}.

\end{document}